\documentclass[12pt]{article}

\title{\textbf{Irreducible Complexity in\\ Pure Mathematics}}

\author{Gregory Chaitin\thanks{IBM Watson Research Center, P. O. Box 218,
Yorktown Heights, NY 10598, USA, chaitin@us.ibm.com}}

\date{}

\begin{document}

\maketitle

\begin{abstract}
By using ideas on complexity and randomness originally suggested by 
the mathematician-philosopher Gottfried Leibniz in 1686,
the modern theory of algorithmic information is able to show
that there can never be a ``theory of everything'' for all of mathematics.
\end{abstract}

\section*{}

In 1956 an article by Ernest Nagel and James R. Newman entitled ``G\"odel's proof''
was published in \emph{Scientific American,}
and in 1958 an expanded version of this article was published as a book with the same title.
This is a wonderful book, and it's still in print. 
    
At the time of its original publication, I was a child, not even a teenager,
and I was obsessed by this little book. I remember the thrill of discovering the newly
published \emph{G\"odel's Proof} in the New York Public Library. 
I used to carry it around with me and try to explain
it to other children.  
    
Why was I so fascinated? Because Kurt G\"odel uses mathematics to show that mathematics itself
has limitations. How can this be? How could reason have limits?
In fact, G\"odel refutes the position of David Hilbert, who about a century
ago declared that there was a theory of everything (TOE) for math, a finite set of principles from
which one could mindlessly deduce all mathematical truths by merely tediously
following the rules of symbolic mathematical logic. 
Such a theory is called a \emph{formal axiomatic mathematical theory.}
    
My attempt to understand G\"odel's proof took over my life, and almost half a century
later, I have just finished writing a little book of my own. It's my own 
version of Nagel and Newman's \emph{G\"odel's Proof,} in which everything 
is done completely differently.
The only thing the two books have in common is their common goal of providing an auto-critique
of mathematical methods and the fact that they are both small books.
My book is currently on my personal website and will soon be published.
    
Why did I have to completely rewrite Nagel and Newman?  Because their exposition and G\"odel's
original 1931 proof are both based on the two self-referential paradoxes: 
``This statement is false'' and ``This statement is unprovable.''
[\emph{See the appendix on ``G\"odel's Proof.''}]
    
My approach is completely different. It's based on measuring information
and on showing that some mathematical facts have no redundancy and cannot be compressed into
any mathematical theory because these facts are too complicated, in fact, infinitely complex.
This new approach suggests that what G\"odel originally discovered was just the
tip of the iceberg and that the problem is much bigger than most people think.
And, amazingly enough, I have recently discovered that these ideas can be traced back
to G. W. Leibniz in the late 17th century.
    
In this essay I'm going to try to outline some of these ideas.  Please bear with me, since
it is not really possible to present all the mathematical details in a short essay. Instead
I'll just try to convey the general flavor of what is going on, using vague but hopefully
suggestive analogies and much hand-waving. If you want to know more, you should read my new book,
\emph{Meta Math!,}
which is as non-technical as I could make it, but takes the time to do things right.
And there are a lot of other books on this subject, by me
and by other people.
    
My story begins in 1686 with Leibniz's philosophical essay 
\emph{Discours de m\'etaphysique,} in English, \emph{Discourse on
Metaphysics,} where Leibniz discusses how we can distinguish between facts that follow
a law, and lawless, irregular, chaotic facts.
How can we do this?
    
Leibniz's idea is very simple and very profound.
It's in section VI of the \emph{Discours.}
It's the observation that the concept of \emph{law} becomes 
vacuous if arbitrarily high mathematical complexity is permitted, for then there is always a law.
Conversely, if the law \textbf{has to be} extremely complicated, 
then the data is irregular, lawless,
random, unstructured, patternless, and also incompressible and irreducible.
A theory has to be simpler than the data that it explains, otherwise it doesn't explain anything.

\section*{Algorithmic Information}

We live in an age of digital information: software, DNA, DVD's, digital cameras, etc.
And the basic idea of algorithmic information theory is to look at the \emph{size}
in bits of computer software, the minimum possible size, without caring at all about
the \emph{speed} of this software.
    
Here is the basic insight, the basic model.
It's a software view of science: A scientific \emph{theory} is a computer program that computes
our \emph{observations,} the experimental data. 
And these are our two fundamental principles, originally due to William of Occam and to Leibniz:
The simplest theory is best (Occam's razor). This means that the
smallest program that calculates the observations is the best theory.
Furthermore,
if a theory is the same size in bits as the data it explains, then it's worthless, because
there is always such a theory (Leibniz).
A theory is a compression of the data, comprehension is compression.
And you compress things into computer programs, into concise algorithmic descriptions.
The simpler the theory, the better you understand something.
A very complex theory means something is wrong.
    
Then you define the \emph{complexity} or, more precisely, the ``algorithmic
information content,'' of a set of facts (or any other digital object) to
be the size in bits of the smallest program for calculating them (or it), no matter how slow.
And a stream of bits is \emph{irreducible} (also called ``algorithmically random'')
if its complexity is equal to its size, in other words, if there is no good theory for it,
no program for calculating it that is smaller than it is.  In other words, such bit streams
are incompressible, they have no redundancy; 
the best thing to do is to transmit them directly. You gain nothing
from doing what DVD players and digital television attempt to do, 
which is to transmit instructions for reassembling
the picture frames rather than the pictures themselves.
    
How close did Leibniz come to these modern ideas? Well, very close.  He had all the key
elements, he just never connected them.
He knew that everything can be represented with binary information, he built one of the
first calculating machines, he appreciated the power of computation, 
and he discussed 
complexity and randomness.  All the key components of the modern theory of algorithmic
information.
    
And if he had put all of this together, he might have been able to realize that one of
the key pillars of his philosophical position can be severely questioned.
For Leibniz, like Spinoza and Descartes, was a \emph{rationalist,}
he believed in the power of reason, he followed that continental philosophical tradition
rather than the British school of empiricism, which is more impressed with facts than
with thoughts and theories.
    
And one of the key elements of Leibniz's philosophy is what he called ``the principle of
sufficient reason,'' that everything happens for a reason. In other words, if something
is true, it must be true for a reason.
Hard to believe sometimes, in the confusion and chaos of daily life, 
in the contingent ebb and flow of human history.
Ah, said Leibniz, even if we cannot always see a reason, perhaps because the chain of
reasoning is long and subtle, God can see the reason, it's there! 
The universe is rational.
In fact, this is actually an idea that comes from the ancient Greeks. 
    
Mathematicians certainly believe in reason and in Leibniz's principle of sufficient
reason, because they always try to prove everything.  No matter how much computational or
experimental evidence there is for something, like the celebrated Riemann hypothesis 
or the $\mathbf{P} \neq \mathbf{NP}$ 
conjecture of theoretical computer science, mathematicians demand a proof,
nothing less will satisfy them!
    
And here is where the concept of algorithmic information can make its surprising contribution
to epistemology, to the philosophical discussion of the origins and limits of knowledge.
What if we can find mathematical facts that are true for no reason, where would that
leave our philosophy, what would it do to us?
   
In fact we can find irreducible mathematical facts, an infinity of them.
Later in this essay I'll show you how. I'll exhibit
an infinite irreducible stream of yes/no, true/false mathematical facts.
These facts turn out to be not only computationally
irreducible, they are even logically irreducible.  In essence, the only way to prove such 
mathematical facts, is to directly assume them as new mathematical axioms,
without using reasoning at all.
    
The concept of
``axiom'' is closely related to the idea of logical irreducibility.  Axioms are mathematical facts
that we take as self-evident and do not attempt to prove from simpler principles.
All formal mathematical theories start with axioms, and then deduce the consequences
of these axioms, which are called its theorems. That is how Euclid did things in Alexandria
two millennia ago, and his treatise on geometry is the classical model for mathematical
exposition, one of the all-time best-sellers of the Western World,
even though Euclid's \emph{Elements}
is no longer taught in high-school the way that it was when I was a child.
    
In ancient Greece 
if you wanted to convince your fellow citizens to vote with you on some issue,
you had to reason with them.
Which I guess is
how they came up with the idea that in math you had to prove things rather than just
discover them experimentally, which is all that it appears that
previous cultures in Mesopotamia and Egypt did.
And using reason has certainly been an extremely fruitful approach, leading to modern mathematics and
mathematical physics and all that goes with it, including, eventually, the technology for
building that extremely logical and mathematical machine, the computer.
    
So am I saying that in a way all this crashes and burns?! 
Yes, in a sense I am. My counter-example illustrating the limited power of logic and reason,
my source of an infinite stream of irreducible/unprovable 
mathematical facts, is the number that I call
$\Omega$, capital omega, the last letter of the Greek alphabet. O-mega as opposed to o-micron, 
which means big-oh, not little-oh. Omega  
is typically used in the expression ``from alpha to omega,'' from the beginning
to the end, or to denote something inaccessible or final, such as the omega point or the omega man.
    
So now let me tell you about $\Omega$ and why it provides us with an infinite stream of
irreducible bits and unprovable facts.

\section*{The Number $\Omega$}

In a famous paper published in 1936 in the \emph{Proceedings of the London Mathematical
Society,} Alan Turing began the computer age by presenting a mathematical model
of a simple general-purpose programmable digital computer.
He then showed that there was something that no algorithm could accomplish and
no mathematical theory could ever decide, namely whether or not a self-contained computer
program will eventually halt.
This is Turing's famous \emph{halting problem.}
    
Of course, by running the program you can eventually discover that it halts.
The problem, and it is an extremely fundamental one, is to decide when to give up.
There is no general solution, even though a great many special cases can in fact
be solved.
And Turing showed in 1936 that there will never be a general solution; it's impossible.
[\emph{For a modern proof, see the appendix ``Why is Turing's halting problem unsolvable?''}]
     
Well, that's the first step on the road to $\Omega$.
    
The next step, which is often a fruitful approach, is to forget about individual cases
of this problem, and consider the statistical ensemble.
In other words, let's ask if a program chosen at random ever halts.
The probability of having that happen is my $\Omega$ number.
     
How precisely do you go about picking a program at random? Well, what you do is
that every time your computer asks for the next bit of the program, you just flip a coin.
So $\Omega$ is just the probability that a machine will eventually do something.
Not a big deal!
    
For those of you who want to know more
about how to define $\Omega$,
let me 
remark
that $\Omega$ is actually the sum of 1/(two raised to the power which is the size
in bits of each program that halts).  
In other words, each $N$-bit program that halts contributes precisely $1/2^N$ to the sum
that defines $\Omega$.
Because if each bit of a program is chosen using an independent toss of a fair coin,
head, 1, tails, 0, then the probability of getting any individual $N$-bit program is precisely
$1/2^N$.
And we want to know the total probability of all the programs that halt.
    
In other words, each $N$-bit program that halts adds a 1 to the $N$th bit in the binary expansion
of $\Omega$.  Do this for all programs that halt, and you'd get the precise value of $\Omega$.
[\emph{For an example, see the appendix ``How to calculate the halting probability $\Omega$.''}]
This may make it sound like you can calculate $\Omega$ with arbitrarily great accuracy, just
as if it were $\sqrt{2}$ or the number $\pi$. However this is actually impossible, because
in order to do it you'd have to solve Turing's halting problem. 
In fact, in a sense the halting probability $\Omega$ is a \emph{maximally uncomputable} number.
     
So $\Omega$ cannot be computed, but it can be defined mathematically as a specific number.
For this to work
it is crucial that no extension of a valid program is a valid program, which
guarantees that the sum 
over all programs that halt defining
$\Omega$ converges to a value between zero and one instead of diverging to
infinity.
Such programs are called \emph{self-delimiting,} they're the ones for which a halting
probability can be defined.
Also, I should mention 
that the precise numerical value of $\Omega$ depends on your choice of computer programming
language, but $\Omega$'s surprising properties aren't at all affected by this choice.
    
Okay those are just
details.
    
The next thing that's crucial is for me to tell you how you get an infinite irreducible
stream of bits out of $\Omega$.
     
Well, that's not difficult.  $\Omega$ is a probability, which has to be between zero and one.
Zero would mean no program halts, one would mean all of them do.  $\Omega$ is actually
greater than zero and less than one, because some programs halt and some don't.
Then write, or more precisely, imagine writing $\Omega$ out in binary, in base-two notation.
Just like the number $\pi$ is 3.1415926\ldots\ with an infinite sequence of digits, $\Omega$, if
you knew its precise numerical value, and if you wrote that out in binary, would give you something
like 0.1110100\ldots\ Let's forget about the initial 0 before the decimal point, and just
consider the infinite stream of bits that comes after the decimal point.
    
Well, it turns out that these are our irreducible mathematical facts, this is an
infinite irreducible stream of bits.
The irreducible math facts are whether each individual bit is a 0 or a 1.
There's no way to know that. These turn out to be facts that are true for no reason,
that are unprovable unless you adopt them as new axioms.
    
Why? Well, it's because $\Omega$ squeezes all the redundancy out of individual cases of
Turing's halting problem.  $\Omega$ is the most compact, the most compressed way of giving
you the answer to all cases of the halting problem. 
In fact, knowing the first $N$ bits of $\Omega$ would in
principle enable you to decide whether or not each program up to $N$ bits in size ever halts.
From this it follows that you need an $N$-bit program in order to be able to calculate $N$ bits of $\Omega$.
From this it also follows that you need $N$ bits of axioms in order to be able to determine $N$ bits of
$\Omega$.
In other words, $\Omega$ is logically as well as computationally irreducible.
    
Can you manage to see why? It's certainly not obvious. 
[\emph{See the appendix ``Why is $\Omega$ incompressible?''}]
I remember realizing that this had to be
the case for the first time when I was flying back to the IBM lab in NY 
after a visit to a university in late 1973
or early 1974. It was a beautiful day, 
and at that precise moment I was staring out the window at the Pentagon in
Washington DC. 
    
And the fact that $\Omega$ is irreducible immediately implies that there cannot be a theory of
everything (TOE) for all of mathematics. Why?
    
The basic intuition is that once you squeeze all the redundancy out of anything,
what you are left with is irreducible, and then whether each bit is 0 or 1
is a fact that in a sense is ``true for no reason,'' it's a complete surprise, it's totally unexpected.
So the bits of $\Omega$ are irreducible mathematical truths, they are mathematical
facts that cannot be derived from any principles simpler than they are, and there
are an infinity of them, an infinity of bits of $\Omega$.  
Math therefore has infinite complexity, whereas any individual
TOE only has finite complexity, and cannot capture all the richness of the full world
of mathematical truth, only a finite part of it.
    
But what does it all mean? 
This doesn't mean that proofs are no good.
And I'm certainly not against reason.
Just because some things are irreducible, doesn't mean we should give up using reasoning.
     
But perhaps sometimes you shouldn't try to prove everything.
Sometimes mathematicians should just add new axioms.
That's what you've got to do if you are faced with irreducible facts.
The problem is knowing that they are irreducible!
    
In a way, saying something is irreducible is giving up, saying that it can't ever be
proved.  Mathematicians would rather die than do that,
contrasting sharply with 
their physicist colleagues, who are happy to be pragmatic and to use plausible reasoning
instead of rigorous proof. 
But remember, to avoid an infinite regress, you've got to stop somewhere, you can't
always prove everything from something else. 
And these irreducible principles, which are called axioms, have
always been a part of mathematics, they are not something new.
I'm just saying that there are a lot more of them out there than people suspected.
     
Physicists are certainly willing to add new principles, new scientific laws, in order
to understand new domains of experience.
Should mathematicians do the same thing that physicists do and add new axioms because
they are pragmatically justified, even though they are not at all ``self-evident?''
This raises what I think is an extremely interesting question:

\section*{Is Math like Physics?}

I don't want to get into a long discussion about this, because it's very, very
controversial, but in fact there are
a few of us who do think that math is like physics, sort of.
    
One of the people who does, or did, is Imre Lakatos, who left Hungary in 1956
and worked on philosophy in England. 
There Lakatos came up with a great word,
\emph{quasi-empirical.}
That's the word that I use to describe my philosophical position on all of this:
I think that math is quasi-empirical.
In other words, I feel that math is different from physics,
but perhaps not as different as most people think.
    
In science you compress your experimental observations into scientific laws,
and in math you compress your computational experiments into mathematical axioms.
What counts is the compression, that's what understanding is, 
it's the fact that you're putting just a few ideas in,
and getting a lot more out, you're using them to explain a rich and diverse
set of physical or mathematical experiences.
    
If Hilbert had been right, math would be a closed system, with no new ideas
coming in.  There would be a static closed finished theory of everything for all of math,
and this would be like a dictatorship.   In fact, for math to progress you actually need 
new ideas and plenty of room
for creativity, it does not suffice to mindlessly and mechanically grind away deducing all the
possible consequences of a fixed number of well-known basic principles.
And I much prefer an open system, I don't like rigid, totalitarian, authoritarian ways of thinking.
    
To end on a personal note, I've lived in both the world of math and
the world of physics, and I never thought there was such a big difference between
these two fields. It's a matter of degree, of emphasis, not an absolute
difference. After all, math and physics co-evolved.
Mathematicians shouldn't isolate themselves.
They shouldn't cut themselves off from rich sources of new ideas.

\footnotesize
\section*{FURTHER READING}

\begin{enumerate}
\item
\textbf{G\"odel's Proof.} (Revised Edition)
E. Nagel, J. R. Newman and D. R. Hofstadter.
New York University Press, 2001.
\item    
For a chapter on Leibniz, see:
\textbf{Men of Mathematics.} E. T. Bell. Simon and Schuster, 1986.
\item
For more on a quasi-empirical view of math, see:
\textbf{New Directions in the Philosophy of Mathematics.}
T. Tymoczko. Princeton University Press, 1998.
\item
For more on $\Omega$, see Chaitin's \emph{Meta Math!}\ e-book:\\
\textbf{http://www.cs.auckland.ac.nz/CDMTCS/chaitin/omega.html}\\
To be published by Pantheon Books in New York in 2005.
\item
For how far you can get in math without bothering with proofs, see:
\textbf{Mathematics by Experiment.}
J. Borwein and D. Bailey. A. K. Peters, 2004.
\textbf{Experimentation in Mathematics.}
J. Borwein, D. Bailey and R. Girgensohn. A. K. Peters, 2004.
\end{enumerate}

\section*{APPENDIXES}

\section*{What is G\"odel's proof?}

Let's start with the paradox of the liar: ``This statement is false.''
This statement is true if and only if it's false,
and therefore is neither true nor false.
    
Now let's consider ``This statement is unprovable.'' 
If it is provable, then we are proving a falsehood, which is extremely unpleasant
and is generally assumed to be impossible.
The only alternative left is that this statement is unprovable.
Therefore it's in fact both true and unprovable, and mathematics is incomplete,
because some truths are unprovable.
   
G\"odel's proof constructs self-referential statements indirectly, using their G\"odel numbers,
which are a way to talk about statements and whether they can be proved by talking about
the numerical properties of very large integers 
that represent mathematical assertions and their proofs.
    
And G\"odel's proof actually shows that what is incomplete is not \emph{mathematics,}
but individual formal axiomatic mathematical theories that pretend to be theories of everything,
but in fact fail to prove the true numerical statement ``This statement is unprovable.''
    
The key question left unanswered by G\"odel: Is this an isolated phenomenon, or are there
many important mathematical truths that are unprovable?

\section*{Why is Turing's halting problem unsolvable?}

A key step in showing that incompleteness is natural and pervasive was taken by Alan Turing 
in 1936 when he demonstrated that there can be no general procedure to decide if a self-contained
computer program will eventually halt.
   
Why not? Well, let's assume the opposite of what we want to prove, namely that there is 
in fact precisely such a general procedure $H$, and from this we shall derive a contradiction.
This is what's called a \emph{reductio ad absurdum} proof.
    
So assuming the existence of $H$, we can construct the following program $P$ 
that uses $H$ as a subroutine.
The program $P$ knows its own size in bits $N$ (there is certainly room in $P$ for it to contain the number
$N$) and then using $H$, which $P$ contains, $P$ takes a look at all programs
up to a hundred times $N$ bits in size in order to see which halt and which don't.  
Then $P$ runs all the ones that
halt in order to determine the output that they produce.
This will be precisely the set of all digital objects with complexity up to a hundred times $N$.
Finally our program $P$ outputs the smallest positive integer not in this set,
that is, with complexity greater than a hundred times $N$, and then $P$ itself halts.
    
So $P$ halts, $P$'s size is precisely $N$ bits, and 
$P$'s output is the first positive integer having 
complexity larger than a hundred times $N$, that is,
that cannot be produced by a program whose size is less than or equal to a hundred times $N$ bits.
However $P$ has just produced this highly complex integer as its output,
and $P$ is much too small to be able to do this, because its size is only $N$ bits, 
which is much less than a hundred times $N$.  Contradiction!
Therefore a general procedure $H$ for deciding whether or not programs ever halt cannot exist,
for if it did then we could actually construct this paradoxical program $P$ using $H$.
    
Finally Turing points out that if there were a TOE 
that always enables you to prove that an individual program halts or to prove that it never does,
whichever is the case,
then by systematically running through all possible proofs you could eventually decide
whether individual programs ever halt. In other words, we could use this TOE to construct $H$, 
which we have just shown cannot exist.
Therefore there is no TOE for the halting problem.

\section*{How to calculate the halting probability $\Omega$}

Let's look at an example.  Suppose that the computer that we are dealing with has only
three programs that halt, and they are the bit strings 0001, 000001 and 000011.
These programs are, respectively, 4, 6 and 6 bits in size.
And the probability of getting each of them by chance is precisely $1/2^4$,
$1/2^6$ and $1/2^6$, since each particular bit has probability $1/2$.
So the value of the halting probability $\Omega$ for this particular computer is precisely
\[
   \Omega 
   \;\;=\;\; \frac{1}{2^4} + \frac{1}{2^6} + \frac{1}{2^6} 
   \;\;=\;\; \frac{1}{2^4} + \frac{2}{2^6} 
   \;\;=\;\; \frac{1}{2^4} + \frac{1}{2^5} 
   \;\;=\;\; .000 110
\]
We're adding a 1 bit into the 4th, 6th and 6th bits of this $\Omega$.
Note that we've gotten a carry out of the 6th bit position after the decimal point.
Adding a 1 bit twice into the 6th bit position 
yields a 0 bit there and a 1 bit in the 5th bit position.
The final value of $\Omega$ is the probability of getting one of these three programs by chance.
     
Here's an embarrassing question.  What if we have a computer for which the two 1-bit
programs 0 and 1 \textbf{both} halt?
Then we would have
\[
   \Omega \;\;=\;\; \frac{1}{2} + \frac{1}{2} \;\;=\;\; 1
\]
This is embarrassing because I said that $\Omega$ had to be greater than zero and less than one.
What is going on here?  Well, remember that I said that programs must be self-delimiting?
This means that the computer decides by itself when to stop asking for more bits, each
of which is determined by an independent toss of a fair coin. 
    
The fact that programs are self-delimiting implies that no extension of a valid
program is a valid program.
Therefore if 0 and 1 were both programs that halt, then no other program could ever halt,
and our computer would only be good for running two programs. That's not much of a computer!
The halting probability $\Omega$ is only of interest, it only has surprising properties, when 
it is the halting probability of what is called a ``universal'' computer. That's
a general-purpose computer that can run all possible algorithms, of which there are 
infinitely many.
    
This is an example of how the fact that no extension of a valid program is a valid program
keeps $\Omega$ from being greater than 1.

\section*{Why is $\Omega$ irreducible/incompressible?}

Well, let's assume that we are given the first $N$ bits of the base-two numerical value of $\Omega$. 
So we know $\Omega$ with accuracy one over two to the $N$, in other words, with $N$-bit accuracy.
Our strategy is to show that this would tell us a lot about Turing's halting problem, which
in fact we already know cannot be solved.
    
More precisely, 
if we can use the first $N$ bits of $\Omega$
to solve the halting problem for all programs up to $N$ bits in size,
this will show that the first $N$ bits of $\Omega$ cannot be produced by a program substantially 
less than $N$ bits
in size nor deduced using a formal axiomatic math theory with substantially less than $N$ bits of 
axioms.
    
Why does knowing $N$ bits of $\Omega$ enable us to solve the halting problem for all programs up to
$N$ bits in size?
We can do this by performing a computation in stages, one for each $K = 1, 2, 3, \ldots\;$
At stage $K$ run every program up to $K$ bits in size for $K$ seconds.
Then compute \textbf{a lower bound} $\Omega_K$ on the halting probability $\Omega$ 
based on all the programs that halt that you discover at stage $K$.
This will give you an infinite list of lower bounds $\Omega_K$ on the actual value of $\Omega$.
    
These lower bounds on $\Omega$ will eventually get closer and closer to the actual value of $\Omega$.
And as soon as the first $N$ bits are correct, you know that you've encountered every program
up to $N$ bits in size that will ever halt,
for otherwise the lower bounds on $\Omega$ would then become larger than 
the true value of $\Omega$, which is impossible.
     
It should be mentioned that the stage $K$ at which the first $N$ bits of $\Omega$ are correct
grows immensely quickly, in fact, faster than any computable function of $N$.

\section*{New mathematical axioms}

There are many well-known examples of \emph{controversial} axioms, including:
\begin{itemize}
\item
the parallel postulate in Euclidean geometry, 
\item
the law of the excluded middle in logic,
\item
the axiom of choice in set theory.
\end{itemize}
But what are examples of potential \emph{new} axioms?
    
Elsewhere the author of this essay has proposed that some variant of
the Riemann hypothesis might be pragmatically justified as a new axiom
due to its diverse applications in number theory.
     
But there are much better examples of new axioms. Here are two that have 
emerged by group consensus in the past few years:
\begin{itemize}
\item
the axiom of projective determinacy in set theory,
\item
the $\mathbf{P} \neq \mathbf{NP}$ conjecture regarding
time complexity.
\end{itemize}
The latter example is particularly interesting as it has many important
applications, as I'll now explain.
     
Most theoretical computer scientists are concerned with time complexity---the
time required to compute things---not with program-size (or informational) complexity as in this
essay. And at the present time this community believes that many important
problems require an amount of time that necessarily grows exponentially,
even though no one can prove it.
     
An example of this presumed exponential growth is determining whether or not
a logical expression involving Boolean connectives is a tautology, that is,
true for all possible assignments of truth values to its variables.
This can be done via truth tables in an exponential amount of time by
looking at all possible combinations of truth values.
That this exponential time growth is necessarily the case is a consequence of
the $\mathbf{P} \neq \mathbf{NP}$ hypothesis, currently believed and freely used by
almost all people working in the field of time complexity.

\section*{Experimental mathematics}

This is the idea of discovering new mathematical results by looking at many
examples using a computer. While this is not as persuasive as a \emph{short}
proof---but it may be more convincing than a long and extremely complicated proof---for 
some purposes it is quite sufficient.
     
Such calculations, checking many diverse examples of a mathematical problem, are
usually done using a symbolic programming language such as Mathematica or Maple,
or using a numerical programming language such as MATLAB.
    
In the past this approach was defended with great vigor by George P\'olya
and by Imre Lakatos, believers in heuristic reasoning and in the quasi-empirical
nature of mathematics.
    
In this generation, experimental mathematics has been promoted most forcefully
by Jon Borwein and David Bailey, authors of a two-volume treatise on the power of
experimental methods.  Another eminent practitioner of experimental math is
my IBM colleague Benoit Mandelbrot, of fractal fame.  This methodology is also
practiced and justified in Stephen Wolfram's \emph{A New Kind of Science.}
    
Contemporary mathematicians are also fortunate to have a journal, \emph{Experimental
Mathematics,} where they can publish their numerical observations and conjectures.
    
Extensive computer calculations can be extremely persuasive, but do they render proof
unnecessary?! Yes and no.  In fact, they provide a \emph{different} kind of evidence.
In important situations, I would argue that \emph{both} kinds of evidence are required,
as proofs may be flawed, and conversely computer searches may have the bad luck to stop
just before encountering a counter-example.

\end{document}